\documentclass[12pt]{amsart}
\usepackage{amstext, amssymb, amscd}

\newfont{\TimeS}{ptmr scaled\magstephalf}
\newfont{\IT}{ptmri scaled\magstephalf}
\newfont{\BF}{ptmb scaled\magstephalf}
\newfont{\BBFF}{ptmb scaled\magstep4}
\newfont{\EM}{ptmrre scaled\magstephalf}

\newcommand{\al}{\ensuremath{\alpha}}
\newcommand{\ga}{\ensuremath{\gamma}}
\newcommand{\be}{\ensuremath{\beta}}

\newcommand{\de}{\ensuremath{\delta}}

\newcommand{\lrw}{\ensuremath{\longrightarrow}} 
\newcommand{\rw}{\ensuremath{\rightarrow}}  
\newcommand{\hrw}{\ensuremath{\hookrightarrow}}  
\newcommand{\op}{\ensuremath{\oplus}}
\newcommand{\BO}{\ensuremath{\bigoplus}}

\newcommand{\ul}{\ensuremath{\underline}}
\newcommand{\wt}{\ensuremath{\widetilde}}

\newcommand{\bs}{\boldsymbol}
\newcommand{\ds}{\displaystyle}
\newcommand{\PP}{\ensuremath{\mathbb P}}

\newcommand{\KK}{\ensuremath{\rm{\bf K}}}
\newcommand{\FF}{\ensuremath{\rm{\bf F}}}

\newcommand{\tT}{\ensuremath{\mathtt T}}

\newcommand{\Hr}{\ensuremath{\mathrm H}}

\newcommand{\codim}{\ensuremath{\mathrm{codim}}}
\newcommand{\grade}{\ensuremath{\mathrm{grade}}}

\newcommand{\xdeg}{\ensuremath{\deg_{\underline x}}}

\newcommand{\kx}{\ensuremath{k[x_1,x_2,\ldots ,x_r]}}

\title[Inverse Gr\"obner basis problem in codimension two]{Inverse Gr\"obner basis problem in codimension two}
\author[A.~Taylor]{{\vspace{0.3cm} Amelia~Taylor${}^*$}}
\address{Department of Mathematics, University of Kansas, Lawrence, Kansas, 66045 U.S.A.}
\email{{\tt ataylor@math.ukans.edu}}

\newtheorem{theorem}{Theorem}[section]
\newtheorem{prop}[theorem]{Proposition}
\newtheorem{lemma}[theorem]{Lemma}

\newtheorem{defin}[theorem]{Definition}
\newtheorem{ex}[theorem]{Example}

\newtheorem{cor}[theorem]{Corollary}

\numberwithin{equation}{section}

\begin{document}

\begin{abstract}
Generic linkage is used to compute a prime ideal such that the radical of the initial ideal of the prime ideal is equal to the radical of a given codimension two monomial ideal that has a Cohen-Macaulay quotient ring.
\end{abstract}
\maketitle
\section*{Introduction}

The \emph{inverse Gr\"obner basis problem} is to find the ideals that have a given monomial ideal as its initial ideal.  We consider the problem of finding when the given monomial ideal is the initial ideal of a prime ideal.   Sturmfels and Kalkbrenner \cite[Theorem 1]{KS} give necessary conditions for a 
square-free monomial ideal to be the radical of the initial ideal of a 
prime ideal.  In particular they prove that the radical of the initial ideal of a prime ideal is equidimensional and connected in codimension one or equivalently, the initial complex of a prime ideal is pure and strongly connected.  They ask if these conditions are sufficient.  

Dalbec \cite[Theorem 2]{D} proved in 1998 that if $I$ is an ideal generated by all the degree $d$ square-free monomials in $n$ variables, then there exists a prime ideal $P$ such that the radical of the initial ideal of $P$ is $I$.  The ideals he considers are all square-free, the generators have the same degree and their quotient is Cohen-Macaulay.  We prove, for codimension two monomial ideals, that the quotient ring being Cohen-Macaulay is sufficient to construct a prime ideal such that the radical of the initial ideal of the prime ideal is equal to the radical of the original ideal.
\begin{theorem}\label{thm:main}Let $R=\kx$ be a polynomial ring over a field $k$ and $I$ an ideal of $R$.  Let $>$ be a monomial order that respects total 
degree.  Assume $I$ is monomial of codimension 2 and $R/I$ is Cohen-Macaulay.  
Then there exists an extension field $K$ of $k$ and a prime ideal $P$ 
contained in the polynomial ring $S=K[x_1,x_2,\ldots ,x_r]$ such that $\sqrt{in(P)}=\sqrt{I}S$.
\end{theorem} 
If $I$ is a square-free monomial ideal and has minimal generators of the 
same degree, then the proof of Theorem~\ref{thm:main} actually gives 
$in(P)=IS$.  If $I$ is only square-free then $\sqrt{in(P)}=IS$.  We would like to have $P$ in $\kx$.  For this we need to specialize indeterminates.  The process of specializing does not necessarily preserve the property of the ideal being prime nor the structure of the initial ideal.  A Bertini theorem~\cite{Fle} can be used to preserve the prime property.  Equations needed to preserve the Gr\"obner basis, in this setting,  are given in~\cite{AT}.  However, the two are incompatible and obtaining $P$ in $\kx$ remains open.

The proof of our main result uses two important ingredients.  
First, generic linkage is 
\vfill
\footnoterule
${}^*${\footnotesize This research was partially conducted by the author for the Clay Mathematics Institute}
\newpage

\noindent
our tool for constructing a prime ideal.  Second, we give a Gr\"obner basis 
for the ideal of maximal minors for a particular class of non-generic matrices.  Gr\"obner bases for ideals of minors of generic matrices are known, however, if the matrix is not generic, 
finding a Gr\"obner basis for the ideal of maximal minors is, in general, difficult.  

We use generic linkage to construct our prime ideal, so in Section 1 we collect the relevant definitions and propositions needed from the theory of generic linkage.  We prove a key theorem on Gr\"obner bases for the ideal of maximal minors for certain non-generic matrices in Section 2.  In Section 3 we prove the main theorem.  Last, in Section 4 we give several examples, including examples that illustrate the construction, explore the necessity of the Cohen-Macaulay hypothesis and the necessity of the radicals. 

Before we proceed to section 1 we include the basic definitions and notation from Gr\"obner basis theory that we will use.
A \emph{monomial order} $\geq$ on a polynomial ring $R=\kx$ over a field $k$ is a total order on the monomials in $R$ such that $m\geq 1$ for each monomial $m$ in $R$ and if $m_1,m_2,n$ are monomials in $R$ with $m_1\geq m_2$ then $nm_1\geq nm_2$.  
A monomial order on a polynomial ring in several variables generalizes the notion of degree for a polynomial ring in one variable.  The \emph{initial term} of an element $f\in R$, denoted 
$in(f)$, is the largest term (including coefficient) of $f$ with respect to a fixed monomial order.  We use $lm(f)$ to denote the largest monomial of $f$ when we do not want to include the coefficient.  Given an ideal $I$ of $R$, the \emph{initial ideal} is defined to be $\langle\{in(f):f\in I\}\rangle$, and is denoted $in(I)$.  
It should be noted that different monomial orders may yield different initial 
ideals, so whenever an initial ideal is referred to, it is assumed a monomial 
order has been fixed.  A \emph{Gr\"obner basis} is a subset $\{g_1, \ldots, g_n\}$ of $I$ such that 
$in(I)=\langle in(g_1),\ldots,in(g_n)\rangle$.


\section{Generic Linkage}

Two varieties $X$ and $Y$ in $\PP^n$ with no common components 
are \emph{linked} if $X\cup Y$ is a complete intersection \cite{PS}.  
Algebraically, two ideals $I$ and $J$ in a local Cohen-Macaulay ring are linked
 if there exists a regular sequence $\al_1,\ldots,\al_s=\ul{\al}$ contained in 
the intersection $I\cap J$ such that $\langle\ul\al\rangle:I=J$ and 
$\langle\ul\al\rangle:J=I$ \cite[Definition 2.1]{HU}.  In \cite{HU}, Huneke 
and Ulrich define a \emph{generic link} of $I$ and prove, under some 
hypotheses, that it is a prime ideal.  
\begin{defin}\label{def:link}{\rm \cite[Definition 2.3]{HU} Let $R$ be a Gorenstein Noetherian 
ring and $I$ an unmixed ideal of $R$ 
of grade $g$.  For $I=R$ we take $g>0$ arbitrary and finite, although the convention in this case is $\grade(I)=\infty$.  Fix a generating set 
$f_1,\ldots, f_m$ of $I$.  A generic link $L(\ul f)$ of $I$ is defined as 
follows:  Let $Y_{ij} (1\leq i\leq g, 1\leq j\leq m)$ be $g\cdot m$ variables 
and set $S=R[Y_{ij}]$ and $\al_i=\sum_{j=1}^m Y_{ij}f_j$, $1\leq i\leq g$.  We 
set $L(\ul f)=\langle\al_1,\ldots,\al_g\rangle : IS$, and call $(S,L(\ul f))$ 
a generic link to $I$.}
\end{defin}

Hochster proved that that $\al_1,\ldots,\al_g$ is a maximal regular sequence in $IS$ \cite{Hoc}.   Therefore if $R$ is Gorenstein and $I$ is unmixed, then $(S,L(\ul f))$ is linked to $I$ \cite{HU}.  Hochster~\cite{Hoc} also gives an equivalence relation on pairs $(R,I)$, where $I$ is an ideal of the ring $R$.  The pairs $(R_1,I_1)$ and $(R_2,I_2)$ are equivalent if there exist integers $r,s$ and indeterminates $Y_1,\ldots, Y_r$ over $R_1$  and $Z_1,\ldots, Z_s$ over 
$R_2$ such that $(R_1/I_1)[Y_1,\ldots, Y_r]$ and $(R_2/I_2)[Z_1,\ldots, Z_s]$ are isomorphic.  Huneke and Ulrich \cite{HU} prove that for $R$ a Gorenstein Noetherian ring and $I$ an unmixed ideal of $R$, if $\ul f$ and $\ul h$ are two generating sets of $I$ and $(Q_1,L(\ul f))$ and $(Q_2,L(\ul h))$ two generic links for $I$, then $(Q_1,L(\ul f))$ is equivalent to $(Q_2,L(\ul h))$.
This alleviates the dependence of Definition~\ref{def:link} on the generating set of $I$ and allows us to use the notation $L(I)$ for a generic link 
of $I$ when the ring is understood.  We use $L_n(I)$ to denote the $n^{th}$ generic 
link of $I$, defined inductively to be $L_1(L_{n-1}(I))$.  The next proposition is the main property we need from the theory of generic linkage.  
\begin{prop}{\rm \cite[Proposition 2.6]{HU}}\label{prop:prime} Let $R$ be a Gorenstein local
domain and let $I$ be an unmixed ideal of $R$ which is generically a complete intersection. 
Let $(S,L(I))$ be a generic link to $I$.  Then $L(I)$ is a prime ideal.
\end{prop}

For any square-free monomial ideal $I\subseteq R=\kx$, $R/I$ is reduced and the primary decomposition looks like $I=P_1\cap\cdots\cap P_s$ where the $P_i$ are generated by subsets of the variables $\{x_1,\ldots,x_r\}$.  Hence any square-free monomial ideal $I\subseteq R$ is generically a complete intersection and $L(I)$ is a prime ideal by Proposition~\ref{prop:prime}.  We use this in the proof of Theorem~\ref{thm:main}.

In the context of the main theorem, $R$ is a polynomial ring and $I$ is a 
homogeneous ideal so $I$ has a graded minimal free resolution.  
Under these assumptions we can construct a free resolution of the generic
 link of $I$.  Set $S$ to be the ring for a generic link of $I$ and $\KK=\KK(\ul{\al};S)$ to be the Koszul resolution of 
$S/\langle\ul{\al}\rangle$.  Let $\FF$ be the minimal free resolution of $S/I$.  A free 
resolution of $S/L(I)$ is the mapping cone of the dual of the map 
$u:\KK\rw \FF$ induced by $S/\langle\ul{\al}\rangle\rw S/I$ \cite[Proposition 2.6]{PS}.  
This resolution has length $\grade(I)+1$, but the last differential in the mapping cone splits.   Taking the mapping cone of the dual of $u$ modulo the subcomplex $S\rw S$ gives a resolution of length equal to the grade of $I$.  

If $R/I$ is Cohen-Macaulay, $I=\langle f_1,\ldots, f_m\rangle$ and $\codim (I)=2$ then, by the Hilbert-Burch theorem \cite[Theorem 1.4.17]{BH}, the resolution $\FF$ has the form
$$
0\lrw R^{m-1}\stackrel{A}\lrw R^m\stackrel{B}\lrw R\lrw 0
$$
\noindent where $B=\begin{bmatrix}(-1)f_1 & (-1)^2f_2 & \ldots & (-1)^mf_m\end{bmatrix}$ and the $(m-1)\times (m-1)$-minors of $A$ generate the ideal $I$, that is $I_{m-1}(A)=I$.  Take $\{Y_{ij}\}_{1\leq j\leq m, 1\leq i\leq 2}$ and form the linear combinations $\al_1=\sum_{j=1}^{m}Y_{1j}f_j$ and $\al_2=\sum_{j=1}^{m}Y_{2j}f_j$.  

In this particular case, after we mod out by $S\rw S$, the resolution for $S/\langle \al_1,\al_2\rangle:IS$ is
\begin{equation}\label{eq:l1res}
0\lrw S^{m}\stackrel{A'}\lrw S^{m+1}\stackrel{B'}\lrw S\lrw 0
\end{equation}
where 
\begin{equation}
        A'=\begin{bmatrix}   & A^{\tT} &   \\ 
                      Y_{11} &\cdots   & Y_{1m}\\
                      Y_{21} &\cdots   & Y_{2m}\end{bmatrix}\quad 
\begin{matrix}\text{ and }B'=\begin{bmatrix}(-1)\de_1 & (-1)^2\de_2 &\ldots & (-1)^{m+1}\de_{m+1}\end{bmatrix}\\ \text{ where }\de_1,\ldots\de_{m+1} \text{ are the signed}\\
\text{ maximal minors of }A'.\end{matrix}
\end{equation}
Since $\al_1, \al_2$ is a regular sequence in $IS$ we know that $\langle \al_1,\al_2\rangle:IS$ has grade at least 2.  Hence its quotient is Cohen-Macaulay and the Hilbert-Burch theorem implies the maximal minors of $A'$ generate $\langle \al_1,\al_2\rangle:IS$, the first generic link of $I$. 

Repeating this process, the second generic link is generated by the maximal minors of the matrix
\begin{equation}\label{eq:thematrix}
A''=\begin{bmatrix}   &       & A^{'\tT} \\ 
               Z_{11} &\cdots & Z_{1m-1} & Y_{1m+1} & Y_{2m+1}\\
               Z_{21} &\cdots & Z_{2m-1} & Y_{1m+2} & Y_{2m+2}\end{bmatrix}=
\begin{bmatrix}   &  A    & Y_{11}  & Y_{21}\\
                  &       & \vdots  & \vdots\\
                  &       & Y_{1m}  & Y_{2m}\\
           Z_{11} &\cdots & Y_{1m+1}  & Y_{2m+1}\\
           Z_{21} &\cdots & Y_{1m+2}  & Y_{2m+2}\end{bmatrix}.\end{equation}
The indeterminates, $Z_{11},\ldots, Z_{2m-1},Y_{1m+1}, Y_{1m+2}, Y_{2m+1}, Y_{2m+2}$, used in forming the second generic link are labelled this way because it is useful in later sections.

We utilize the second generic link and its resolution in the proof of Theorem~\ref{thm:main}.  There are two 
reasons to expect the second generic link to be better than the first generic link for our purposes.  First, $I$ is a link 
of $L_1(I)$ and $L_2(I)$ is a generic link of $L_1(I)$, so there is a specialization from the second generic link of $I$ to $I$~\cite[Proposition 2.13]{HU}.

Second, the structure of the maximal minors of $A''$ (\ref{eq:thematrix}) is better.  We give an example.  Let $f_i$ denote the signed minor of $A$ when the $i^{th}$ row is removed , for $1\leq i\leq m$, and let $\de_i$ denote the same signed minor for $A''$.  Therefore,
\begin{equation}\label{eq:L2}\de_i=f_i(Y_{1m+2}Y_{2m+1}-Y_{1m+1}Y_{2m+2})+\be_i\quad 
1\leq i\leq m.\end{equation}
where each term of $\be_i$, $1\leq i\leq m$ has higher degree in the new variables.  For example, the following are the first three minors of a presentation matrix for the second generic link of $I=\langle ac,ad,bd\rangle$.  

\begin{equation}\begin{split}\label{ex:good}
\de_1 =& {\bf ac(Y_{15}Y_{24}-Y_{14}Y_{25})}-aY_{15}Y_{22}Z_{12}+aY_{12}Y_{25}Z_{12}+aY_{14}Y_{22}Z_{22}\\
       &{}-aY_{12}Y_{24}Z_{22}-bY_{15}Y_{23}Z_{12}+bY_{13}Y_{25}Z_{12}+bY_{14}Y_{23}Z_{22}-bY_{13}Y_{24}Z_{22}\\
       &{}-cY_{15}Y_{23}Z_{11}+cY_{13}Y_{25}Z_{11}+cY_{14}Y_{23}Z_{21}-cY_{13}Y_{24}Z_{21}+Y_{13}Y_{22}Z_{12}Z_{21}\\
       &{}-Y_{12}Y_{23}Z_{12}Z_{21}-Y_{13}Y_{22}Z_{11}Z_{22}+Y_{12}Y_{23}Z_{11}Z_{22}\\
\de_2 =& {\bf ad(Y_{15}Y_{24}-Y_{14}Y_{25})}+aY_{15}Y_{21}Z_{12}-aY_{11}Y_{25}Z_{12}-aY_{14}Y_{21}Z_{22}\\
       &{}+aY_{11}Y_{24}Z_{22}-dY_{15}Y_{23}Z_{11}+dY_{13}Y_{25}Z_{11}+dY_{14}Y_{23}Z_{21}-dY_{13}Y_{24}Z_{21}\\
       &{}-Y_{13}Y_{21}Z_{12}Z_{21}+Y_{11}Y_{23}Z_{12}Z_{21}+Y_{13}Y_{21}Z_{11}Z_{22}-Y_{11}Y_{23}Z_{11}Z_{22} \\
\de_3 =& {\bf bd(Y_{15}Y_{24}-Y_{14}Y_{25})}+bY_{15}Y_{21}Z_{12}-bY_{11}Y_{25}Z_{12}-bY_{14}Y_{21}Z_{22}\\
       &{}+bY_{11}Y_{24}Z_{22}+cY_{15}Y_{21}Z_{11}-cY_{11}Y_{25}Z_{11}-cY_{14}Y_{21}Z_{21}+cY_{11}Y_{24}Z_{21}\\
       &{}+dY_{15}Y_{22}Z_{11}-dY_{12}Y_{25}Z_{11}-dY_{14}Y_{22}Z_{21}+dY_{12}Y_{24}Z_{21}+Y_{12}Y_{21}Z_{12}Z_{21}\\
       &{}-Y_{11}Y_{22}Z_{12}Z_{21}-Y_{12}Y_{21}Z_{11}Z_{22}+Y_{11}Y_{22}Z_{11}Z_{22} 
\end{split}\end{equation}
The non-boldface terms are the terms we call $\be_i$ in Equation (\ref{eq:L2}) and the boldface terms correspond to the remaining in Equation (\ref{eq:L2}).  

The minors $\de_{m+1}$ and $\de_{m+2}$ do not have exactly the form given in Equation (\ref{eq:L2}), but as can be seen below each has terms of the form $f_iM$ where $M$ is a degree two monomial in the new variables.  
\begin{equation*}\begin{split}
\de_4 =& {\bf ac(Y_{15}Y_{21}-Y_{11}Y_{25})+ad(Y_{15}Y_{22}-Y_{12}Y_{25})+bd(Y_{15}Y_{23}-Y_{13}Y_{25})}\\
       &{}-aY_{12}Y_{21}Z_{22}+aY_{11}Y_{22}Z_{22}-bY_{13}Y_{21}Z_{22}+bY_{11}Y_{23}Z_{22}-cY_{13}Y_{21}Z_{21}\\
&{}+cY_{11}Y_{23}Z_{21}-dY_{13}Y_{22}Z_{21}+dY_{12}Y_{23}Z_{21} \\
\de_5 =& {\bf ac(Y_{14}Y_{21}-Y_{11}Y_{24})+ad(Y_{14}Y_{22}-Y_{12}Y_{24})+bd(Y_{14}Y_{23}-Y_{13}Y_{24})}\\
       &{}-aY_{12}Y_{21}Z_{12}+aY_{11}Y_{22}Z_{12}-bY_{13}Y_{21}Z_{12}+bY_{11}Y_{23}Z_{12}-cY_{13}Y_{21}Z_{11}\\
       &{}+cY_{11}Y_{23}Z_{11}-dY_{13}Y_{22}Z_{11}+dY_{12}Y_{23}Z_{11}       
\end{split}\end{equation*}
In this example, where the degrees of the generators of $I$ have the same degree, the initial term of each $\de_i$, $1\leq i\leq m+2$, using the inverse block order with respect to the added variables (see the beginning of Section 2 for a definition of this order), is of the form $f_iM$ where $M$ is a degree two monomial in the new variables.  We will use this structure in the proof of Theorem~\ref{thm:main} as well as the proof of Theorem~\ref{thm:gbasis}.
 
In contrast, the minors of a presentation matrix for the first generic link of $I$ do not have this form.
The following are the maximal minors of a presentation matrix for the first generic link of $I=\langle ac,ad,bd\rangle$.  
\begin{equation}\begin{split}\label{eq:L1}
\de_1 =& -cY_{12}Y_{21}+cY_{11}Y_{22}+dY_{13}Y_{22}-dY_{12}Y_{23}\\
\de_2 =& aY_{13}Y_{22}-aY_{12}Y_{23}+bY_{13}Y_{21}-bY_{11}Y_{23}\\
\de_3 =& acY_{21}+adY_{22}+bdY_{23}\\
\de_4 =& acY_{11}+adY_{12}+bdY_{13}
\end{split}\end{equation}
The first and second minors are linear in $a,b,c,d$.  These linear terms will appear in any Gr\"obner basis computed from this generating set.  The linear terms are not in $I$ and all of this suggests that the second generic link, as opposed 
to the first, is the prime ideal we want to work with.

\section{Gr\"obner bases for ideals of maximal minors}
In this section we prove, under certain conditions, that the maximal minors of the presentation matrix of $L_2(I)$ given in equation (\ref{eq:L2}) are a Gr\"obner basis for $L_2(I)$.  
First, we need to fix some notation.  
From now on assume that the monomial order on $R=\kx$ respects total degree.  
The \emph{inverse block order} \cite[$\S 8$]{KW} is a useful monomial 
order when adding additional variables to $R$.  Let $<_R$ denote the monomial 
order on $R$.  Let $<_T$ denote the monomial order on $T=k[y_1,\ldots,y_s]$.  
Set $S=k[x_1,\ldots,x_r,y_1,\ldots,y_s]$ and $Y=\{y_1,\ldots,y_s\}$.  Let 
$s_1,s_2\in R$ and $t_1,t_2\in T$ be monomials.  The \emph{inverse block 
order on $S$ with respect to $Y$} is defined as follows:  if either $s_1<s_2$, 
or $s_1=s_2$ and $t_1<t_2$, then $s_1t_1<s_2t_2$.  Given a monomial $st$ in 
$S$, where $s\in R$ and $t\in T$, define $\xdeg(st)=\deg s$.  

For $I\subseteq R$ homogeneous use $\Hr(R/I,n)=\dim_k(R/I)_n$ to denote the Hilbert 
function for $R/I$ and $\Hr_{R/I}(t)=\sum_{n=0}^{\infty}\Hr(R/I,n)t^n$ to denote the Hilbert series for $R/I$.  We use the following standard fact.  If $\FF$ is a minimal graded free resolution of $R/I$, then the Hilbert series is the alternating sum of the Hilbert series of the modules in the resolution, so $\Hr_{R/I}(t)=\sum_{i=0}^m(-1)^i\Hr_{F_i}(t)$ where $m=pd_R(R/I)$.  

We also use the following standard facts in the setup and proof of Theorem~\ref{thm:gbasis}.  Let $f_1,\ldots, f_m$ be a generating set for a homogeneous ideal $I$ of $R$.  If $I$ is 
codimension two and $R/I$ is Cohen-Macaulay then the Hilbert-Burch theorem along with \cite{PS} implies that a minimal graded free resolution of $I$ has the following form

\begin{equation}0\lrw \BO_{j=1}^{m-1}R(-b_j)\stackrel{\phi}\lrw \BO_{i=1}^{m}R(-a_i)\stackrel{A}\lrw R\lrw R/I\lrw 0.\end{equation}
Denote the $ij^{th}$ entry of $\phi$ by $\phi_{ij}$  Furthermore~\cite{PS},   
\begin{gather}
a_j=\deg(f_j),\quad \text{for } 1\leq j\leq m,\label{eq:ps1}\\ 
\deg(\phi_{ij})=b_j-a_i\text{, for the non-zero entries of }\phi_{ij},\label{eq:ps2}\\
A=\begin{bmatrix} (-1)f_1 & (-1)^2f_2 & \cdots & (-1)^mf_m\end{bmatrix},\label{eq:ps3}\\
I=I_{m-1}(\phi)\label{eq:ps4}.
\end{gather}
If $\deg(f_1)=\deg(f_2)=\cdots=\deg(f_m)=d$, then Equation (3.3) implies $\deg(\phi_{ij})=b_j-d$ 
for the non-zero entries of $\phi$.  Hence all of the non-zero entries in the $j^{th}$ 
column of $\phi$ have degree $b_j-d$, $1\leq j\leq m-1$.  And vice versa, if all of the non-zero entries in each column of $\phi$ have the same degree the maximal minors must be homogeneous and of the same degree.

Before proceeding to Theorem~\ref{thm:gbasis} we indicate a way to simplify computations.  
The need for this can be seen in the minors of the presentation matrix of the second generic
link of $I=\langle ac,ad,bd\rangle$ (\ref{ex:good}).  The first generic link may involve simpler computations, but yields polynomials that appear unhelpful (see Equation~\ref{eq:L1}) so we simplify in a different way.  Let $\phi$ be a presentation matrix for $I$.  The following matrix is the matrix that is useful for simplifying computations.     
\begin{equation}\label{matrix:thm}
  \Phi= \left( \begin{array}{cccc}       &         &         & Y_{1}  \\ 
                                         &  \phi   &         & \vdots  \\
                                         &         &         & Y_{m}  \\
                                  Z_{1}  & \cdots  & Z_{m-1} & Y_{m+1}   
\end{array} \right).\end{equation} 
The minors of matrix~\ref{matrix:thm} for the ideal $I=\langle ac,ad,bd\rangle$ are included to illustrate how much simpler they are while maintaining the structure of the minors of the second generic link.  The boldface terms are the initial terms using the inverse block order and the remaining terms are what we call $\be_i$ in the proof of Theorem~\ref{thm:gbasis}.  Compare these equations to Equation (\ref{ex:good}).
\begin{equation}\begin{split}\label{ex:thm}
\de_1 =& {\bf acY_{14}}-aY_{12}Z_{12}-bY_{13}Z_{12}-cY_{13}Z_{11} \\
\de_2 =& {\bf adY_{14}}+aY_{11}Z_{12}-dY_{13}Z_{11}\\
\de_4 =& {\bf bdY_{14}}+bY_{11}Z_{12}+cY_{11}Z_{11}+dY_{12}Z_{11}\\  
\de_5 =& {\bf acY_{11}}+adY_{12}+bdY_{13}
\end{split}\end{equation}
Theorem~\ref{thm:gbasis} is about the minors of this matrix.  In Corollary~\ref{cor:L2} we add a second row and column of new variables to prove that under the given conditions, the maximal minors of the presentation matrix for the second generic link, form a Gr\"obner basis for the second generic link.
The following definition will aid in stating the next theorem and corollary.  

\begin{defin}\label{def:propA}          
{\rm Let $R=\kx$ and $I=\langle f_1,\ldots,f_m\rangle\subseteq R$. \\
We say the pair $(R,\underline{f})$ has \emph{Property $A$} if\\ 
(1) $\codim (I)=2$.\\
(2) $R/I$ Cohen-Macaulay.\\
(3) $I$ homogeneous.\\ 
(4) $deg(f_1)=\cdots=deg(f_m)$.\\ 
(5) $\{f_1,\ldots, f_m\}$ is both a Gr\"obner basis and a minimal generating set for $I$.}
\end{defin}
The first two conditions in this definition are the assumptions for the main theorem.  Conditions 3 and 5 are satisfied by any monomial ideal and help in the induction that follows.  Condition 4 is seemingly strong, however we reduce the main theorem to this case using Corollary~\ref{cor:simpledet}.  This condition allows us to say that the minors of $\Phi$ are homogeneous so we can use graded free resolutions.  In Section 4 (see page 20) we give an example of what happens, besides no-longer having graded resolutions, if we do not make this assumption.
\begin{theorem}\label{thm:gbasis}Let $R=\kx$ and 
$I=\langle f_1,\ldots,f_m\rangle\subseteq R$.  Assume $(R,\underline{f})$
 has Property $A$.  Let $\phi$ be a graded free presentation matrix for $I$ and $e_j$ the degree 
of the non-zero entries in the $j^{th}$ column of $\phi$.  Let $Y=\{Y_1,\ldots,Y_{m+1}\}$, 
$Z=\{Z_1,\ldots,Z_{m-1}\}$ and $S=R[Y,Z]$.  Give $S$ the inverse block order with respect to 
$\{Y,Z\}$.  Set $\deg(Y_i)=1$ for $1\leq i\leq m+1$ and $\deg(Z_j)=e_j$ for $1\leq j\leq m-1$.  Set 
$$\Phi= \left( \begin{array}{cccc}       &         &         & Y_{1}  \\ 
                                         &  \phi   &         & \vdots  \\
                                         &         &         & Y_{m}  \\
                                  Z_{1}  & \cdots  & Z_{m-1} & Y_{m+1}   
\end{array} \right).$$ 
Then the maximal minors of $\Phi$ form a Gr\"obner basis for the ideal $I_{m}(\Phi)$.   
\end{theorem}

\proof For $1\leq i\leq m+1$, let $\de_i$ denote the signed minor of $\Phi$ when 
the $i^{th}$ row is removed.  The generators, $f_1,\ldots, f_m$ of $I$ are a Gr\"obner basis by assumption, so $in(f_1),\ldots, in(f_m)$ generate $in(I)$.  
We prove that the Hilbert series for $S/I_m(\Phi)$ and $S/\langle in(\de_1),\ldots, in(\de_{m+1})\rangle$ are equal.  Therefore $in(I_m(\Phi))=\langle in(\de_1),\ldots, in(\de_{m+1})\rangle$ 
and hence the maximal minors of $\Phi$ form a Gr\"obner basis for $I_m(\Phi)$. 

We can order $f_1,\ldots,f_m$ such that $in(f_1)\geq in(f_2)\geq \ldots \geq in(f_m)$ and 
assume $f_i$ is the signed $(m-1)\times(m-1)$ minor of $\phi$ when the $i^{th}$ row is removed.  Let $d=\deg(f_i)$ for $1\leq i\leq m$.  

The non-zero entries in each column of $\Phi$ have the same degree, so by the 
remarks before the theorem, $\de_i$, $1\leq i\leq m+1$, is homogeneous and has 
the form        
\begin{equation}\label{eq:delta}\begin{split}
\de_1     &=f_1Y_{m+1}+\be_1\\
\vdots    &\\
\de_m     &=f_mY_{m+1}+\be_n\\
\de_{m+1} &=f_1Y_1+f_2Y_2+\cdots +f_mY_m.
\end{split}\end{equation}
Each $\be_i$ is at least degree two in the new variables so $\xdeg(\be_i)<d$, $1\leq i\leq m$.  Also 
$\codim(I_m(\Phi))\leq 2$ \cite[Theorem 3]{EN}.
Suppose $\codim(I_m(\Phi))\leq 1$.  Localize at a codimension one prime ideal 
$P$ containing $I_m(\Phi)$.  Since $I_{m-1}(\phi)$ is codimension two, at 
least one $(m-1)\times(m-1)$ minor of $\phi$ is invertible in $S_P$.  Using 
row and column operations, rewrite $\Phi$ as  
$$\Phi'= \left( \begin{array}{cccc}  1    &   0     & \cdots   & 0   \\ 
                                    0    &   1     & \cdots   & 0 \\
                                  \vdots & \vdots  & \vdots   & \vdots \\
                                    0    &   0     &   1      & 0   \\
                                    0    &   0     &   0      & Y_{m}'  \\
                                    0    & \cdots  &   0      & Y_{m+1}'   
\end{array} \right).$$
The terms $Y_{m}'=Y_m+\be$ and $Y_{m+1}'=Y_{m+1}+\ga$ are such that $\be$ and $\ga$ are polynomials not involving $Y_m$ and $Y_{m+1}$ respectively.  Then $Y_{m}'$ and $Y_{m+1}'$ are a regular sequence and $I_m(\Phi)S_P=I_m(\Phi')S_P=\langle Y_m',Y_{m+1}'\rangle S_P$ is a codimension two 
ideal in $S_P$.  This is a contradiction and hence $\codim (I_m(\Phi))=2$.  Since $S$ is 
a regular ring, $\grade(I_m(\Phi))=\codim(I_m(\Phi))=2$.  Therefore, the Hilbert-Burch Theorem implies 
$$0\lrw S^{m}\stackrel{\Phi}\lrw S^{m+1}\stackrel{A}\lrw S\lrw S/J\lrw 0$$
is a minimal free resolution for $R/I_m(\Phi)$ and 
$A=\begin{bmatrix}(-1)\de_1 & \cdots & (-1)^{m+1}\de_{m+1}\end{bmatrix}$.  
Each $\de_i$, $1\leq i\leq m+1$, is homogeneous of degree 
$d+1$, by construction, and if $\Phi_{ij}\neq 0$ then $b_j$, as defined in Equation (3.1), is  
\begin{equation*}
b_j=\begin{cases}
     \deg(\Phi_{ij})+(d+1)=e_j+(d+1) & 1\leq j\leq m-1,\\
     d+2 & j=m.
     \end{cases}
\end{equation*}
The following is a graded free resolution of $S/I_m(\Phi)$ with the twists.
$$0\lrw \BO_{j=1}^{m-1}S(-d-e_j-1)\op S(-d-2)\stackrel{\Phi}\lrw S(-d-1)^{m+1}\stackrel{A}\lrw S\lrw S/J\lrw 0.$$
Therefore 
$$\Hr_{S/I_m(\Phi)}(t)={\dfrac{1-(m+1)t^{d+1}+t^{d+2}+\sum_{j=1}^{m-1}t^{d+e_j+1}}{(1-t)^N(1-t^{e_1})(1-t^{e_2})\cdots(1-t^{e_{m-1}})}},$$
where $N=r+m+1$.

Each $\be_i$, $1\leq i\leq m+1$, in Equation (\ref{eq:delta}) has the property that $\xdeg(\be_i)<d$.  Also, $R$ has an order that respects total degree and $S$ has the inverse block order, therefore $\langle in(\de_1),\ldots,in(\de_{n+1})\rangle=\langle in(f_1)Y_{m+1},\ldots ,in(f_m)Y_{m+1},in(f_1)Y_1\rangle$.  Let $K$ denote this ideal.  
The following is a standard exact sequence.
  \begin{equation}\label{eq:seq}
0\lrw {\dfrac{S}{(K:Y_{m+1})}}(-1)\lrw {\dfrac{S}{K}}\lrw {\dfrac{S}{(K,Y_{m+1})}}\lrw 0
\end{equation}
Since $\langle K,Y_{m+1}\rangle=\langle Y_{m+1},in(f_1)Y_1\rangle$ and $(K:Y_{m+1})=in(I)$, Equation (\ref{eq:seq}) gives
\begin{equation}\label{eq:sum}
\Hr_{S/K}(t)=\Hr_{S/(Y_{m+1},in(f_1)Y_1)}(t)+\Hr_{S/in(I)(-1)}(t).
\end{equation}

The monomials $Y_{m+1}$ and $in(f_1)Y_1$ consist of distinct variables and therefore form a regular sequence.  Hence
\begin{equation}\label{eq:regular}\Hr_{S/(Y_{m+1},in(f_1)Y_1)}(t)={\dfrac{1-t-t^{d+1}+t^{d+2}}{(1-t)^N(1-t^{e_1})(1-t^{e_2})\cdots(1-t^{e_{m-1}})}}.\end{equation}
Viewing $I$ as an ideal in $S$, the 
Hilbert-Burch Theorem and Equation (\ref{eq:ps2}) provide the following graded free resolution of $S/I$.  
\begin{equation*}
0\lrw \BO_{j=1}^{m-1}S(-d-e_j)\stackrel{\phi}\lrw S(-d)^m\lrw S\lrw 0.
\end{equation*}
Therefore, 
\begin{equation}\label{eq:original}{\dfrac{1-mt^d+\sum_{j=1}^{m-1}t^{d+e_j}}{(1-t)^{N}(1-t^{e_1})(1-t^{e_2})\cdots(1-t^{e_{m-1}})}}=\Hr_{S/I}(t)=\Hr_{S/in(I)}(t),\end{equation}
where the last equality is standard \cite[Theorem 15.26]{Eis}.  To shift this series by 1, multiply by $t$.  Combining Equations (\ref{eq:sum}), (\ref{eq:regular}) and (\ref{eq:original}) 

\begin{gather}
\Hr_{S/K}(t)={\dfrac{1-t-t^{d+1}+t^{d+2}+t-mt^{d+1}+\sum_{j=1}^{m-1}t^{d+e_j+1}}{(1-t)^N(1-t^{e_1})(1-t^{e_2})\cdots(1-t^{e_{m-1}})}}=\\
={\dfrac{1-(m+1)t^{d+1}+t^{d+2}+\sum_{j=1}^{m-1}t^{d+e_j+1}}{(1-t)^N(1-t^{e_1})(1-t^{e_2})\cdots(1-t^{e_{m-1}})}}=H_{S/I_{m}(\Phi)}(t). 
\end{gather} \qed
\medskip

\begin{cor}\label{cor:L2}
Assume $R$, $I$, $S$, $\phi$, $\Phi$ and $\{\de_i\}_{i=1}^{m+1}$ are as in Theorem~\ref{thm:gbasis}.  Then \\
{\rm (1)} $(S,\underline{\de})$ has Property $A$\\
{\rm (2)} The generators for $L_2(I)$ from the presentation matrix (\ref{eq:thematrix}) form a Gr\"obner basis for $L_2(I)$.
\end{cor}

\proof (1):  In the proof of the theorem we established that $I_m(\Phi)$ is 
codimension two and $S/I_m(\Phi)$ is Cohen-Macaulay.  The construction of 
$\Phi$ implies $\de_i$ is homogeneous of degree $d+1$, for $1\leq i\leq m+1$.  
No entry in $\phi$ is a unit by assumption, so the same is true of $\Phi$ by 
construction.  Hence $\{\de_i\}_{i=1}^{m+1}$ form a minimal generating set for 
$I_m(\Phi)$.  Theorem~\ref{thm:gbasis} implies $\{\de_i\}_{i=1}^{m+1}$ is a 
Gr\"obner basis for $I_m(\Phi)$.

(2):  By (1), we can apply Theorem~\ref{thm:gbasis} to $(S,\underline{\de})$.  
The matrix that arises in this process is the same as the presentation matrix 
for $L_2(I)$ given in Section 1 (see~\ref{eq:thematrix}) and hence the maximal 
minors form a Gr\"obner basis for $L_2(I)$. \qed

\medskip


\section{Sufficient conditions in Codimension 2}

For a monomial ideal $I=\langle f_1,\ldots,f_m\rangle$ in $R=\kx$ inductively define a \emph{polarization} of $I$ as follows.  Let $\al_{j}$ denote the exponent of $x_1$ in $f_j$ for $1\leq j\leq m$.  Write each $f_j$, $1\leq j\leq m$, as $x_1^{\al_{j}}m_j$ where $x_1$ does not divide $m_j$.  Set $\ds{\al=\max_{1\leq j\leq m}\{\al_{j}\}}$ and $Y_{1},\ldots,Y_{\al-1}$ to be $\al-1$ new indeterminates.  Set
\begin{equation} 
P(f_j)=
\begin{cases}
x_1Y_{1}Y_{2}\cdots Y_{\al_{j}-1}m_j &\text{if } \al_{j}\geq 2,\\
f_j                                  &\text{if } \al_{j}=0,1.
\end{cases}\end{equation}
For each $f_j$, $1\leq j\leq m$, repeat this process for each $x_i$, $1\leq i\leq r$, and call the resulting monomial the \emph{polarization of $f_j$}.  A \emph{polarization of $I$} is the ideal generated by the polarizations of $f_1,\ldots,f_m$.  Let $P(\underline{f})$ denote the polarization of $I$, formed from the generating set $\underline{f}=f_1,\ldots,f_m$.  

The polarization is a square free monomial ideal by construction.  Let ${\bf Y}$ denote the indeterminates needed to form the polarization of $I$ and set 
$${\bf Y}-{\bf x}=\{Y_i-x_j\vert Y_i \text{ replaces } x_j \text{ in the polarization}\}.$$
Then $R[{\bf Y}]/\langle P(\underline{f}),{\bf Y}-{\bf x} \rangle\simeq R[{\bf Y}]/\langle I,{\bf Y}-{\bf x}\rangle.$  The following proposition, which is folklore, uses the polarization of $I$ in its proof.  We use this proposition in our proof of Theorem~\ref{thm:main} to reduce to the case where the monomial ideal $I$ is square-free.
\begin{prop}\label{prop:CM}Let $R=\kx$ and $I$ be a monomial ideal of R.  
If $R/I$ is Cohen-Macaulay then $R/\sqrt{I}$ is also Cohen-Macaulay.
\end{prop}
\proof Set a minimal 
generating set for $I=\langle f_1,\ldots, f_m\rangle$.  First, using induction, we prove that $R[{\bf Y}]/P(\ul{f})$ is Cohen-Macaulay.  Denote the 
degree of $x_1$ in $f_j$ by $\al_j$ for $1\leq j\leq m$.  Reorder the generators of $I$ so that 
$\al_{1},\ldots,\al_{s}>1$, $\al_{s+1}=\cdots=\al_{s+r}=1$ and 
$\al_{s+r+1}=\cdots=\al_{m}=0$.  Write $f_i=x_1^{\al_{i}}h_i$, $1\leq i\leq s+r$, such that $x_1$ does not divide $h_i$.   Use $Y_1$ to denote the first variable used to replace 
$x_1$.  We claim $Y_{1}-x_1$ is a non-zero divisor on the ideal 
$$J=\langle x_1^{\al_{1}-1}Y_{1}h_1,\ldots,x_1^{\al_{s}-1}Y_{1}h_s,x_1h_{s+1},\ldots,x_1h_{s+r},f_{s+r+1},\ldots,f_m\rangle.$$
Suppose $Y_{1}-x_1$ is a zero divisor on $J$.  Thus $Y_1-x_1$ is in some associated prime ideal $P$ of $J$.  There exists $g\notin J$ such that $P=(J:g)$.  Since $J$ is monomial, $P$ is monomial and $g$ can be taken to be monomial.  So there exists a monomial $g\notin J$ such that $g(Y_1-x_1)\in J$.  

Since $J$ is a monomial ideal $Y_1g\in J$ and $x_1g\in J$.  
Since $Y_1$ is a non-zero divisor on $\langle x_1h_{s+1},\ldots,x_1h_{s+r},f_{s+r+1},\ldots,f_m\rangle$, if $Y_1g$ is in that ideal, then $g\in J$ which is a contradiction.
Therefore we can assume $x_1^{\al_i-1}Y_1h_i$ divides $Y_1g$ for some $1\leq i\leq s$.  Write $g=x_1^nY_1^lg'$ where $l,n\geq 0$ and $Y_1$ and $x_1$ do not divide $g'$.  Since $x_1^{\al_i-1}Y_1h_i$ divides $Y_1g$ for some $1\leq i\leq s$, $x_1^{\al_i-1}h_i$ divides $g=x_1^nY_1^lg'$.  Thus, $x_1$ divides $g$ so $n>0$.  Also, $x_1^{\al_i-1}h_i$ divides $x_1^ng'$.  Therefore, $x_1^{\al_i-1}Y_1h_i$ divides $x_1^nY_1g'$.  If $l>0$ this divides $g$ and then $g$ is in J which is a contradiction.  Hence we may assume $n>0$ and $l=0$.  Now we use that $x_1g\in J$.  The monomials $x_1^{\al_i-1}Y_1h_i$ cannot divide $x_1g$ since $Y_1$ does not divide $x_1g$.  
This implies $x_1^ng'\in \langle h_{s+1},\ldots,h_{s+r},f_{s+r+1},\ldots,f_m\rangle$ and $x_1$ is a non-zero divisor on the ideal, so $g'\in \langle h_{s+1},\ldots,h_{s+r},f_{s+r+1},\ldots,f_m\rangle$.  
Since $n>0$, $x_1g'$ divides $g$ and therefore $g\in J$ a contradiction.  Hence $Y_1-x_1$ is a non-zero divisor on $J$.      

Assume $S$ is a graded ring and $m$ is the irrelevant ideal, then $S$ is Cohen-Macaulay
 if and only if $S_m$ is Cohen-Macaulay~\cite{MR}.  
Let $m_1$ denote the irrelevant ideal for $R[Y_1]$.  Since $I$ is homogeneous and $R/I$ is Cohen-Macaulay, the previous two statements combine to imply $R[Y_1]_{m_1}/IR[Y_1]_{m_1}$ is Cohen-Macaulay.  The ring $R[Y_1]_{m_1}/IR[Y_1]_{m_1}$ is local and $Y_1-x_1\in m_1$ is a non-zero divisor, so 
$$R[Y_{1}]_{m_1}/\langle I, Y_{1}-x_1\rangle R[Y_{1}]_{m_1}\simeq R[Y_{1}]_{m_1}/\langle J, Y_{1}-x_1\rangle R[Y_{1}]_{m_1}$$ 
is Cohen-Macaulay.  The ring $R[Y_{1}]_{m_1}/JR[Y_{1}]_{m_1}$ is Cohen-Macaulay since the element $Y_{1}-x_1$ is a non-zero divisor on $J$.  Moreover, this implies 
$R[Y_1]/JR[Y_1]$ is Cohen-Macaulay.  By induction on the variables used to form a polarization of $I$, both $R[{\bf Y}]/P(\underline{f})$ and $R[{\bf Y}]/\langle P(\underline{f}), {\bf Y}-{\bf x}\rangle$ are Cohen-Macaulay.

Let $W$ be the multiplicatively closed set $k[{\bf Y}]\setminus\{0\}$ in $S=R[{\bf Y}]$ and let $K=k({\bf Y})$.  Then the localization of $S/P(\underline{f})$ at $W$ is isomorphic to $K[x_1,\ldots,x_r]/\sqrt{I}$.  So, $S/P(\underline{f})$ 
Cohen-Macaulay implies $K[x_1,\ldots,x_r]/\sqrt{I}$ is Cohen-Macaulay.  The ring   
$$K[x_1,\ldots,x_r]/\sqrt{I}\simeq \kx/\sqrt{I}\otimes_k K$$
is Cohen-Macaulay if and only if $\kx/\sqrt{I}$ is Cohen-Macaulay~\cite[Theorem 2.1.10]{BH}. Hence $\kx/\sqrt{I}$ is Cohen-Macaulay.\qed
\medskip

Recall that Proposition~\ref{prop:prime} states that $L_2(I)$ is a prime ideal assuming $I$ is generically a complete intersection.  In general, monomial ideals with Cohen-Macaulay quotient ring are not generically complete intersections, however, if the ideal is a square-free monomial ideal then it is generically a complete intersection.  Proposition~\ref{prop:CM} allows us to reduce to the square-free case in the proof of Theorem~\ref{thm:main}.

In the statement of Theorem~\ref{thm:main} we do not 
assume the generators of $I$ have the same degree, but this is required by 
Theorem~\ref{thm:gbasis} and we use Theorem~\ref{thm:gbasis} in the proof of Theorem~\ref{thm:main}.  Proposition~\ref{prop:simpledet} and 
Corollary~\ref{cor:simpledet} allow us to reduce to the case where the degrees 
of the generators have the same degree.

For $R=\kx$ and $I=\langle f_1,\ldots, f_m\rangle$ monomial the relations on the generators of $I$ are generated, due to the natural multi-grading, by relations of the form $mf_i-nf_j$ where $m,n\in R$ are monomials .  Hence, 
a presentation matrix of $I$ can be given with exactly two 
non-zero monomial entries in each column.  We call a determinant of an 
$n\times n$ matrix $A=(a_{ij})$ a \emph{simple determinant} if 
at most one term of 
$\det(A)=\ds{\sum_{\sigma\in S_n}sgn(\sigma)a_{1\sigma(1)}\cdots a_{n\sigma(n)}}$ is non-zero.
\begin{prop}\label{prop:simpledet}Let $R=\kx$. Let $\phi$ be a $m\times(m-1)$ matrix with exactly two non-zero monomial entries in each column and assume the maximal minors of $\phi$ are all non-zero.  
Then each minor of $\phi$ is a simple determinant or zero.
\end{prop}
\proof We use induction on the size of $\phi$.  If $m=2$, the minors of $\phi$ are the two monomial entries of $\phi$ and hence are simple determinants.  Assume 
the statement for $m$.  Let $\phi$ be a $(m+1)\times m$ matrix satisfying the hypotheses of the proposition.  Each column has exactly two non-zero entries so there are exactly $2m$ non-zero entries in $\phi$.  Since the maximal minors are all non-zero, every row has at least one non-zero element in it and therefore there must be at least one row with exactly one non-zero entry.  
Choose a row with exactly one non-zero entry and denote that entry $M$.  Set $N$ to be the other non-zero entry in that column.  Reorder the rows and columns so that $\phi$ looks like
\begin{equation*}
\phi =           \begin{bmatrix} M & 0 & 0     & \cdots & 0\\
                                 N &   &       &        &  \\
                                 0 &   & \psi  &        &  \\
                            \vdots &   &       &        &  \\
                                 0 &   &       &        &  \end{bmatrix}
\end{equation*}
where $\psi$ is the $m\times(m-1)$ matrix obtained by removing the first row and column from $\phi$.  
Every column of $\psi$ must have exactly two non-zero entries since otherwise we contradict this fact for $\phi$.  The entries of $\psi$ are monomial since $\psi$ is a submatrix of $\phi$.  Suppose a maximal minor of $\psi$ is zero.  Let $i$ denote the row that was removed to form the maximal minor and let $\de_i$ denote this minor.  The maximal minor of $\phi$ when $i+1^{st}$ row is removed is $M$ times $\de_i$.  Therefore the $i+1^{\rm{st}}$ maximal minor of $\phi$ is $M\de_i=M\cdot 0=0$, which is a contradiction.  Hence $\psi$ satisfies the induction hypotheses and therefore we can assume all of the minors of $\psi$ are simple determinants or zero.  
  Consider the minors of $\phi$.  If a square submatrix of $\phi$ is also a submatrix of $\psi$ then the determinant is simple or zero.  Suppose the submatrix is not contained in $\psi$.  If the submatrix does not contain $M$ or $N$ then, since it is not contained in $\psi$, it must have a row or column of zeros and hence the determinant is zero.   Suppose the submatrix includes $M$ (and may or may not include $N$), if the determinant is expanded along the top row of the submatrix the determinant is $M$ times the determinant of a submatrix contained in $\psi$ and is hence zero or simple.  Last, assume $N$ is in the submatrix, but $M$ is not.  In this case the submatrix of $\phi$ is a submatrix of 
 \begin{equation*}
\begin{bmatrix}                  N &   &       &        &  \\
                                 0 &   & \psi  &        &  \\
                            \vdots &   &       &        &  \\
                                 0 &   &       &        &  \end{bmatrix}.
\end{equation*}
Expanding the determinant along the first column, the determinant is $N$ times the determinant of a submatrix of $\psi$ and hence is zero or simple. \qed
\medskip

This proposition is particularly interesting because given a codimension two monomial ideal $I$ with Cohen-Macaulay quotient ring we can now construct a codimension two monomial ideal $J$ with Cohen-Macaulay quotient ring such that $\sqrt{I}=\sqrt{J}$, the generators of $J$ have the same degree and $J$ is generically a complete intersection.  
\begin{cor}\label{cor:simpledet}Let $R=\kx$ and $I\subseteq R$ be a codimension two monomial ideal such that $R/I$ is Cohen-Macaulay.  Then there exists a monomial ideal $J$ such that $\sqrt{J}=\sqrt{I}$, $R/J$ is Cohen-Macaulay, $J$ is generically a complete intersection and the generators of $J$ are all of the same degree.
\end{cor}

\proof The ring $R/\sqrt{I}$ is Cohen-Macaulay by Proposition~\ref{prop:CM} and if such a $J$ exists for $\sqrt{I}$ then it satisfies the properties for $I$ as well.  Therefore replacing $I$ with $\sqrt{I}$, we can assume $I$ is square-free and generically a complete intersection.

Let $f_1,\ldots,f_m$ be a minimal generating set for $I$.  Let $\phi$ denote 
a presentation matrix for $I$ from the Hilbert-Burch theorem such that there 
are exactly two non-zero entries in each column.  Proposition~\ref{prop:simpledet} allows us to ``homogenize'' $\phi$ in the following way.  Fix a column of $\phi$ and compare the two non-zero entries in that column.  Raise the exponents in the monomial of smaller degree until the two monomials 
have the same degree.  Do this for all of the columns of $\phi$.  Call this 
new matrix $h(\phi)$ (see example~\ref{matrix:homog}).  The non-zero entries 
in each column of $h(\phi)$ are the same degree, by construction, so the 
maximal minors of $h(\phi)$ all have the same degree.  By Proposition~\ref{prop:simpledet}, each minor of both $\phi$ and $h(\phi)$ is a simple determinant 
or zero.  This implies the maximal minors of $h(\phi)$ are monomial.  By 
construction, a minor of $\phi$ is non-zero if and only if the corresponding 
minor of $h(\phi)$ is non-zero.  The construction of $h(\phi)$ and the fact 
that the non-zero minors of each matrix are simple, implies that if $\al$ is a non-zero minor of $\phi$ and $\be$ is the corresponding minor of $h(\phi)$ then $\al$ divides $\be$ and for $N\gg 0$, $\be$ divides $\al^N$.  Hence  
\begin{equation}\label{eq:homogphi}I_{n}(\phi)=\sqrt{I_{n}(h(\phi))}\quad\text{ for }1\leq n\leq m-1.\end{equation}
This implies $\codim (I_{m-1}(h(\phi)))=2$ and hence $R/I_{m-1}(h(\phi))$  is Cohen-Macaulay since $R$ is a regular ring.     

In the context of this corollary $\codim(I_{m-2}(\phi))\geq 3$ if and only if $I$ is generically a complete intersection.  Since $I$ is square-free and monomial, $I$ is generically a complete intersection and hence $\codim(I_{m-2}(\phi))\geq 3$.  By Equation (\ref{eq:homogphi}) $I_{m-2}(\phi)=\sqrt{I_{m-2}(h(\phi))}$.  Therefore $\codim(I_{m-2}(h(\phi)))\geq 3$ and $I_{m-1}(h(\phi))$ is generically a complete intersection.\qed
\medskip

We need one last lemma before proceeding to the proof of the main theorem.

\begin{lemma}\label{lem:ci}Let $R=k[x_1,\ldots,x_r]$ and $I\subseteq R$ an ideal.  Let $(S,L(I))$ be the generic link of $I$.  Let $Y$ denote the variables used to form $L(I)$, so $S=R[Y]$.  If $L(I)\cap k[Y]\neq(0)$ then $I$ is a complete intersection.
\end{lemma}
\proof 
Set $g=\grade(I)$.  Assume there exists $f(Y)\neq 0$ in $L(I)\cap k[Y]$.  Set $W=k[Y]\setminus\{0\}$.  Since $f(Y)\in L(I)$, $L(I)_W$ is the whole ring.  Therefore $\langle \ul{\al}\rangle_W=I_W$, and hence in $k(Y)[\ul{x}]$, $I_W$ is generated by $g$ elements.  Since $k(Y)[\ul{x}]$ is faithfully flat over $k[\ul{x}]$, $I$ had to be generated by $g$ elements to start with and is hence a complete intersection.  \qed
\medskip

We now give the proof of Theorem~\ref{thm:main} and restate the theorem for the reader's convenience.

\newtheorem{them2}{Theorem 0.\hspace{-0.2cm}}
\begin{them2}Let $R=\kx$ be a polynomial ring over a field $k$ and 
$I$ an ideal of $R$.  Let $>$ be a monomial order that respects total degree.  Assume 
$I$ is monomial of codimension 2 and $R/I$ is Cohen-Macaulay.  Then there exists an extension 
field $K$ of $k$ and a prime ideal $P$ contained in the polynomial ring $S=K[x_1,x_2,\ldots ,x_r]$ such that $\sqrt{in(P)}=\sqrt{I}S$.
\end{them2}

\proof By Corollary~\ref{cor:simpledet} there exists a monomial ideal $J$ that is generically a complete intersection, the generators are all the same degree, $\sqrt{J}=\sqrt{I}$ and $R/J$ is Cohen-Macaulay.  If we prove the theorem for $J$ we get a prime ideal $P$ such that $\sqrt{in(P)}=\sqrt{J}S=\sqrt{I}S$ as desired. 

Replace $I$ with $J$ and let $\phi$ denote the presentation matrix for $J$ constructed from a presentation matrix of $I$ as in Corollary~\ref{cor:simpledet}.  Let $m$ denote the number of rows of $\phi$ and $f_i$, $1\leq i\leq m$ be the signed minors of $\phi$.  By the Hilbert-Burch theorem, $f_1,\ldots,f_m$ generate $J$ and by Corollary~\ref{cor:simpledet} they are the same degree.  Denote that degree $d$.  Also, order the generators so that $f_1>f_2>\cdots>f_m$.

By Proposition~\ref{prop:prime} $L_2(I_{m-1}(\phi))$ is a prime ideal.  Set $Y=\{Y_{11},\ldots,Y_{2m+2}\}$, $Z=\{Z_{11},\ldots,Z_{2m-1}\}$ and $Q=R[Y,Z]$ with $\deg(Y_{ij})=1$ and $\deg(Z_{ij})=e_j$.   Give $Q$ the inverse block order with respect to $\{Y,Z\}$ and any order on the new variables.  Set  
\begin{equation}
\Phi=\begin{bmatrix}   &        &         & Y_{11}   & Y_{21}\\
                       &  \phi  &         & \vdots   & \vdots\\
                       &        &         & Y_{1m}   & Y_{2m}\\
                Z_{11} & \cdots & Z_{1m-1}& Y_{1m+1} & Y_{2m+1}\\
                Z_{21} & \cdots & Z_{2m-1}& Y_{1m+2} & Y_{2m+2}. \end{bmatrix}
\end{equation}
The matrix $\Phi$ is a presentation matrix for $L_2(I_{m-1}(\phi))$.  
Let  $\de_i$, $1\leq i\leq m+2$ denote the signed minor of $\Phi$ formed when the $i^{th}$ row 
is removed.  Then,  
\begin{align*}
\de_1     =&f_1(Y_{1m+2}Y_{2m+1}-Y_{1m+1}Y_{2m+2})+\be_1\\
\vdots     &\\
\de_m     =&f_m(Y_{1m+2}Y_{2m+1}-Y_{1m+1}Y_{2m+2})+\be_m\\
\de_{m+1} =&f_1(Y_{11}Y_{2m+2}-Y_{12}Y_{1m+2})+f_2(Y_{12}Y_{2m+2}-Y_{22}Y_{1m+2})+\cdots +\\
           &+f_m(Y_{1m}Y_{2m+2}-Y_{2m}Y_{1m+2})+\be_{m+1}\\
\de_{m+2} =&f_1(Y_{11}Y_{2m+1}-Y_{21}Y_{1m+1})+f_2(Y_{12}Y_{2m+1}-Y_{22}Y_{1m+1})+\cdots +\\
           &+f_m(Y_{1m}Y_{2m+1}-Y_{2m}Y_{1m+1})+\be_{m+2}
\end{align*}  
where $\xdeg(\be_i)<d$, $1\leq i\leq m+2$.  These form a Gr\"obner 
basis for $L_2(I_{m-1}(\phi))$, by Theorem~\ref{thm:gbasis}.  
The images of $\de_1,\ldots,\de_{m+2}$ in 
$S=k(Y,Z)[x_1,\ldots,x_r]$ are a Gr\"obner basis for the image of 
$L_2(I_{m-1}(\phi))$ in $S$ because $Q$ has the inverse block order with respect to $\{Y,Z\}$ \cite[Lemma 8.93]{BW}.  

In the context of this theorem the first generic link is never a complete intersection because it is a grade two ideal and a minimal generating set is the maximal minors of a matrix with $m+1$ rows and $m$ columns, where $m\geq 2$ is the number of generators of $I$.  Therefore $L_2(I_{m-1}(\phi))\cap k[Y,Z]=(0)$ by Lemma~\ref{lem:ci}.  Let $W$ be the multiplicatively closed set $k[Y,Z]\setminus\{0\}$.  Hence $L_2(I_{m-1}(\phi))$ is disjoint from $W$ and the image of $L_2(I_{m-1}(\phi))$ in $Q_W=S$ is a prime ideal.  

The image of $\de_i$, $1\leq i\leq m+2$, in $S$ is $\de_i$ where the part 
of each term that is a monomial in $k[Y,Z]$ is considered to be part of the 
coefficient.  We are using the inverse block order and $\xdeg(\be_i)<\deg(f_i)$ so the initial terms of the $\de_i$ are 
\begin{gather*}
in(\de_1)=f_1Y_{m+2}Y_{m+1}\\
in(\de_2)=f_2Y_{m+2}Y_{m+1}\\
\vdots\\  
in(\de_m)=f_mY_{m+2}Y_{m+1}\\
in(\de_{m+1})=f_1Y_{m+2}Y_1\\
in(\de_{m+2})=f_1Y_{m+1}Y_1.
\end{gather*}
Denote the image of $\de_i$ in $S$ by $\wt{\de_i}$.  Therefore, $lm(\wt{\de_i})=f_i$, for $1\leq i\leq m$ and $lm(\wt{\de_{j}})=f_1$ for $j=m+1,m+2$.  
Hence $in(L_2(I_{m-1}(\phi))S)=\langle f_1,f_2,\ldots,f_n\rangle S$ 
and 
$$IS = (f_1,\ldots,f_m)S=in(L_2(I_{m-1}(\phi))S).$$ \qed
\medskip

\section{Examples}

The ideal $I=\langle bc,bd,acd\rangle$ is a nice test case for many of the theorems and assumptions.  This ideal is one of the simplest examples that illustrates the necessity of each of the steps we have taken.  We will use this ideal for the first and third examples.  The first example illustrates the need for the assumption from Theorem~\ref{thm:gbasis} that the generating set have elements of the same degree, as well as the necessity of one of the radicals.  The second example illustrates the necessity of the other radical.   For the third example we use $I$ to illustrate the entire algorithm for constructing the desired prime ideal.  
The fifth example looks at the necessity of the assumption that the quotient ring be Cohen-Macaulay.  
\begin{ex}\label{ex:1}{\rm
This example illustrates what happens if the assumption in Theorem~\ref{thm:gbasis} that the generators of $I$ have the same degree is dropped.  It also illustrates why we need one of the radicals.  Let $I=\langle bc,bd,acd\rangle$.  Using Macaulay2 \cite{M2} we computed the minors of a presentation matrix for the second generic link for the ideal $I$ as in Theorem~\ref{thm:gbasis}. 
\begin{equation}\begin{split}\label{ex:bad2}
\de_1 & =acY_{13}Z_{11}-{\bf bcY_{14}}+bY_{12}Z_{11}+cY_{13}Z_{12}  \\
\de_2 & =bdY_{14}+bY_{11}Z_{11}-dY_{13}Z_{12}  \\
\de_3 & =acdY_{14}+acY_{11}Z_{11}+cY_{11}Z_{12}+dY_{12}Z_{12}  \\ 
\de_4 & =acdY_{13}+bcY_{11}+bdY_{12}
\end{split}\end{equation}
The elements $\de_1$, $\de_3$ and $\de_4$ are not homogeneous and the bold face term is the one we would like to have as the initial term of $\de_1$.  Set $\deg(Y_{13})=\deg(Z_{11})=1$ and $\deg(Y_{11})=\deg(Y_{12})=\deg(Y_{21})=\deg(Y_{22})=\deg( Y_{14})=\deg(Z_{12})=2$.  The polynomials in (\ref{ex:bad2}) are now quasi-homogeneous, meaning they are homogenous with respect to the weights.  However, the initial term is still not the desired one.  If we now use the reverse lexicographic order with $Y_{14}$ as the largest of the new variables, $\de_1$ now has initial term $-bcY_{14}$.  For any ideal where the generators are not all of the same degree we can weight the new variables so that the second generic link is quasi-homogeneous and use a particular reverse lexicographic order and get a statement like Theorem~\ref{thm:gbasis} for ideals where the generators are not the same degree.  However, in the proof of Theorem~\ref{thm:main}, after constructing the second generic link of $I$, we invert all of the new variables to get an ideal in $K[x_1,\ldots,x_r]$ where $K$ is an extension field of $k$.  After inverting the new variables the leading monomial of the image of $\de_1$ is $ac$ not $bc$ as needed.  This happens regardless of the monomial order and the weights on the new variables. Also, one might suggest sending the new variables to elements in $k$ so that the problematic terms go to zero, however, this process will not preserve the property that the ideal is prime.  Hence $in(P)\neq IS$ if we do not reduce to the case where the generators all have the same degree.  Therefore we raise the degree of $b$ in the presentation matrix using Proposition~\ref{prop:simpledet} and get $\sqrt{in(P)}=IS$. 
}\end{ex}

\begin{ex}\label{ex:radical2}{\rm 
Let $I=\langle x^2,xy,y^2\rangle$, then $I$ is codimension two and the quotient ring $k[x,y]/I$ is Cohen-Macaulay, but I is not generically a complete intersection.  Therefore in order for the second generic link to be a prime ideal we must pass to the radical of $I$, $\langle x,y\rangle$.  Since this ideal is square-free and the generators are of the same degree, the algorithm will yield a prime ideal $P$ with initial ideal $\langle x,y\rangle$ and hence $in(P)=\sqrt{I}$.
}\end{ex}

\begin{ex}{\rm
Using the ideal $I=\langle bc,bd,acd\rangle$ we work through the entire algorithm for constructing the desired prime ideal.  In this example we are able to take the process one step further, as we are able to specialize the new variables to elements in $k$ and verify that the image is indeed a prime ideal.  
In every example we have computed, specializing is possible, however, as we already mentioned, the fact that we can always specialize and preserve both the Gr\"obner basis and the property of being prime is open.  

The following proposition is one way to verify that the ideal we construct is a prime ideal.
\begin{prop}{\rm \cite[Proposition 3.5.6]{V}}\label{prop:wolmer}
Suppose $A=k[{\bf z},{\bf x}]/I$ is a Cohen-Macaulay, equidimensional ring.  Let $B=k[\bf{z}]$ be a Noether normalization of $A$.  
The degree of $A$ over $B$ is the dimension of the vector space
$$l=\dim_k(k[{\bf z},{\bf x}]/(I,{\bf z})).$$
If there exists a subring 
$$B\hrw S=k[{\bf z},U]/\langle f({\bf z},U)\rangle\hrw A,$$
where $f({\bf z},U)$ is an irreducible polynomial of degree $l$, then $A$ is an integral domain.
\end{prop}

First, compute a presentation matrix for $I=\langle bc,bd,acd\rangle\subseteq R=k[a,b,c,d]$.  This ideal is small enough we can find a presentation by hand, or using Macaulay2.  A matrix is
\begin{equation}\phi= \begin{bmatrix} -d & 0 \\ c & -ac \\ 0 & b\end{bmatrix}.
\end{equation}
Since the generators of $I$ are not all of the same degree we use Corollary~\ref{cor:simpledet} and form the following matrix
\begin{equation}\label{matrix:homog}h(\phi)= \begin{bmatrix} -d & 0 \\ c & -ac \\ 0 & b^2\end{bmatrix}\end{equation} 
where $b$ is squared in the second column so that the non-zero entries in each column of $\phi$ now have the same degree.  The maximal minors of $h(\phi)$ are $b^2c,b^2d,acd$, thus verifying that $\sqrt{I_2(h(\phi))}=I_2(\phi)$.  The ideal generated by the maximal minors of $h(\phi)$ is still generically a complete intersection by 
Corollary~\ref{cor:simpledet} and hence the second generic link of $h(\phi)$ is a prime ideal.  Generators for the second generic link are the maximal minors of the following matrix
\begin{equation*}\Phi= 
\begin{bmatrix} -d & 0      & Y_{11} & Y_{21}\\ 
                 c & -ac    & Y_{12} & Y_{22}\\
                 0 & b^2    & Y_{13} & Y_{23}\\
             Z_{11}& Z_{12} & Y_{14} & Y_{24} \\
             Z_{21}& Z_{22} & Y_{15} & Y_{25}\end{bmatrix}.\end{equation*}
Corollary~\ref{cor:L2} implies the maximal minors of this matrix form a Gr\"obner basis for $L_2(I)$.  We can use Macaulay2 to verify this.  
Using Macaulay2 generate a random sequence of numbers from the base field.  For this example we generated a random sequence of rational numbers.  Specialize the new variables $Y_{11},\ldots,Y_{25},Z_{11},\ldots,Z_{22}$ to the following numbers, respectively.  
$$  {\frac {-9} {5}},\quad {\frac{3}{10}},\quad {\frac{-1}{2}},\quad {\frac 7 5},\quad  -1,\quad  -7,\quad  6,\quad  7,\quad  1,\quad  {\frac{-1}{7}} ,\quad{\frac{-1}{4}},\quad {\frac{2}{5}},\quad 1,\quad{\frac{3}{2}}$$
The image of three of the elements in the Gr\"obner basis have $acd$ as their leading monomial.  Therefore, the Gr\"obner basis reduces to include only one of these generators and the Gr\"obner basis simplifies to three elements.  The following three elements form a reduced Gr\"obner basis for the image of $L_2(I)$ after specializing the variables added when the second generic link was formed.
\begin{equation*}\begin{split}
g_1=&\bs{acd}+{\frac{807}{-440}}ac+{\frac{3777}{-1375}}c+{\frac{12123}{-5500}}d-{\frac{51429}{220000}} \\
g_2=&\bs{b^2d}-{\frac{807}{440}}b^2-{\frac{139}{22}}d+{\frac{27993}{-22000}}    \\
g_3=&\bs{b^2c}+{\frac{137}{88}}b^2+{\frac{103}{-22}}ac+{\frac{139}{-22}}c-{\frac{10633}{-22000}} 
\end{split}\end{equation*}

Homogenize with respect to $t$.  Set $S=k[a,b,c,d,t]$ and $G_i=g_i({\frac a t},{\frac b t},{\frac c t})t^3$ for $1\leq i\leq 3$.  
The ideal $\langle g_1,g_2,g_3\rangle$ is a prime ideal if and only if the ideal $\langle G_1,G_2,G_3\rangle$ is a prime ideal.  Let $J=\langle G_1,G_2,G_3\rangle$.  
The degree of $S/J$ is 7.  This can be found either by using Macaulay2 or by computing a Noether Normalization of $S/J$, say $k[\underline{z}]$ and computing $\dim_kS/\langle J,\underline{z}\rangle$, utilizing Proposition~\ref{prop:wolmer}.  
Reorder the variables in $S$ so that $t>a>b>c>d$ and recompute the Gr\"obner basis using an elimination order for $t$.  There is one polynomial in the Gr\"obner basis which is contained in $\langle a,b,c,d\rangle$.  Denote this polynomial $f(a,b,c,d)$.  
A Noether normalization of $S/J$ is $k[a,c,b+d]$.  Hence $$k[a,c,b+d]\hrw k[a,b,c,d]/\langle f(a,b,c,d)\rangle =k[a,b,c,d]/J\cap k[a,b,c,d]\hrw S/J.$$
The polynomial we found has degree 7.  Using Maple we checked that it is irreducible and hence $J$ is a prime ideal.  
}\end{ex}

\begin{ex}{\rm
We consider the necessity of the quotient ring being Cohen-Macaulay.  
Since the two ideals in this example are not perfect we use the definition of generic link to compute the ideals. 

Let $R=k[a,b,c,d,e]$ and $I=(ad,ace,bcd,bce)=(a,b)\cap(a,c)\cap(c,d)\cap(d,e)$.  This ideal is pure and strongly connected which are the necessary conditions given by Kalkbrenner and Sturmfels \cite{KS}.  
Localizing at $P=(a,b,d,e)$, a codimension 2 prime ideal in $R/I$, $R_P/I_P$ is not depth 2 so $R/I$ is not $S_2$.   
The leading terms of a minimal generating set for $L_2(I)$ are given below.
\begin{align*}
b^2cY_{1,3}Z_{1,4}Z_{2,3}        & \qquad \bs{-a^2eY_{1,2}Z_{1,4}Z_{2,3}}\\
abcY_{1,3}Z_{1,4}Z_{2,3}         & \qquad -a^2Y_{1,2}Z_{1,4}Z_{2,3}\\
a^2eY_{1,2}^2Z_{1,4}Z_{2,3}      & \qquad -abcY_{1,2}Z_{1,4}Z_{2,3}\\
\bs{-a^2bY_{1,2}Z_{1,4}Z_{2,3}}   & \qquad a^2Y_{1,2}Z_{1,4}^2Z_{2,3}^2\\
abcY_{1,3}Y_{2,2}Z_{1,3}         & \qquad abcY_{1,3}Y_{2,2}Z_{2,3}
\end{align*}
The portion of each of these monomials that is in $k[a,b,c,d,e]$ is not in $I$.  Moreover, if we carefully check each generator we see that the two elements in the Gr\"obner basis with the boldface monomials as leading terms contain no term whose $a,b,c,d,e$ part is in $I$.  Any Gr\"obner basis will preserve this bad structure.  

The case when $I$ is $S_2$ but not Cohen-Macaulay is both more and less encouraging.  The second generic link in this case does not give a counter example, but we cannot compute it.  Every example we have tried is too computationally complex for the computer we use.  Let $R=k[a,b,c,d,e]$ and $I=(abd,bde,ace,acd,bce)=(a,b)\cap(a,c)\cap(c,d)\cap(d,e)$.  The ideal $I$ is one such example.  
The leading terms for a minimal generating set for the first generic link are not promising, but the first generic link wasn't promising for the Cohen-Macaulay case either.  There is other evidence in a paper by Hochster and Huneke~\cite{HH} that suggests $S_2$ may be the desired necessary condition.
 }\end{ex}

\section*{Acknowledgments}
The results in this paper are from my PhD thesis at the University of Kansas.  I would like to thank my advisor, Craig Huneke, for his help and suggestions.  

%


\end{document}